\documentclass[a4paper,10pt]{article}
\usepackage{calc,amsmath,theorem,amssymb,amsfonts,color}
\usepackage[all,2cell]{xy}
\SelectTips{cm}{10} \UseAllTwocells

\def\labelstyle{\scriptstyle}
\theoremstyle{break}
\newtheorem{theorem}{Theorem}[section]
\newtheorem{definition}[theorem]{Definition}
\newtheorem{proposition}[theorem]{Proposition}
\newtheorem{observation}[theorem]{Observation}
\newtheorem{example}[theorem]{Example}
\newtheorem{corollary}[theorem]{Corollary}
\newenvironment{proof}{\begin{trivlist}\item[]{\em Proof\/}: }
                {\hspace*{\fill}$\blacksquare$\end{trivlist}}

\newenvironment{myenumerate}{%
 \begin{enumerate}
\itemsep 0.33\itemsep plus 0.33\itemsep minus 0.33\itemsep
\parskip .5pt plus 1pt minus .5pt
 }{\end{enumerate}}
\newenvironment{myitemize}{%
 \begin{itemize}
\itemsep 0.33\itemsep plus 0.33\itemsep minus 0.33\itemsep
\parskip .5pt plus 1pt minus .5pt
 }{\end{itemize}}
\renewcommand{\a} {{\mathcal A}}

\renewcommand{\c} {{\mathcal C}}
\newcommand{\e} {{\mathcal E}}
\newcommand{\g} {\ensuremath{{\mathcal G}}}
\renewcommand{\k} {\ensuremath{{\mathcal K}}}

\newcommand{\ba} {{\mathbf A}}
\newcommand{\bb} {{\mathbf B}}
\newcommand{\bh} {\mathbf h}
\newcommand{\Endos} {\text{\sf END}}
\def\endos {\text{\sf End}}
\newcommand{\Ext} {\text{\rm Ext}}
\newcommand{\lalpha} {{\s\alpha}}
\newcommand{\lbeta} {{\s\beta}}
 
\newcommand{\prshv}[1] {{\sf Set}^{\op{#1}}}
\newcommand{\ssetsdel} {\prshv{\sdelta}}
\newcommand{\ssets} {{\sf Sset}} \let\sset=\ssets
\newcommand{\sdelta} {{\Delta}}
\newcommand{\Gr} {{\sf Gp}} \let\gp\Gr \let\Gp\Gr
\newcommand{\iso} {{\sf Iso}}
\newcommand{\Cat} {{\sf Cat}} \let\cat\Cat
\newcommand{\Gpd} {{\sf Gpd}} \let\gpd\Gpd
\newcommand{\dosgpd} {{2\text{-}\Gpd}}
\newcommand{\doscat} {{2\text{-}\Cat}}
\newcommand{\ldoscat} {{2\text{-}\Cat_{\text{\rm lax}}}}
\newcommand{\ldosgpd} {{2\text{-}\Gpd_{\text{\rm lax}}}}
\newlength{\mytmplgth}
\newlength{\mdpt}
\newlength{\mhei}
\newcommand{\mthstrut}[2][0ex]{\settoheight{\mhei}{\ensuremath{#2}}%
                            \settodepth{\mdpt}{\ensuremath{#2}}%
                            \setlength{\mhei}{\mhei+.1ex}%
                            \setlength{\mhei}{\mhei+#1}\rule{0pt}{\mhei}}
\newcommand{\abs}[1]  {{\left|#1\right|}}
\newcommand{\lact}[3][0ex]{{\mthstrut[#1]{#3}}^{#2}{\kern-.3em}{#3}}
\newcommand{\chobj}[3][1]{{\parbox[b][#1\height]{0pt}
                 {\vphantom{\ensuremath{#3}}}}^{#2}{\kern-.3em}{#3}}
\newcommand{\charr}[2] {\til{#1}_{_{\!#2}}} \let\til\widetilde
\newcommand{\id}[1]   {{\mathop{\text{\rm id}}\nolimits}_{#1}}

\newcommand{\op}[1]   {{#1^{\text{\rm op}}}}
\newcommand{\semi}[2] {{\int_#2\kern-.1em#1}}
\newcommand{\LaxFun}  {{\text{\sf LaxFun}}}
\newcommand{\Laxact}  {{\text{\sf Act}}\kern-12pt
                                \widetilde{\phantom{iii}}\kern 2.3pt}
\newcommand{\aut}    {{\text{\sf Aut}}}
\newcommand{\Aut}    {{\text{\sf AUT}}}
\newcommand{\fib}    {{\text{\sf Fib}}}
\newcommand{\bofib}  {\fib_{\text{\rm b.o.}}}

\newcommand{\embed}  {{\text{\bf i}}}
\newcommand{\ner}    {{\text{\sf Ner}}}

\newcommand{\inv}    {{\mathop{\text{\rm inv}}\nolimits}}
\newcommand{\obj}[1] {#1_{\text{obj}}}
\newcommand{\Obj}    {{\text{\sf Obj}}}
\newcommand{\Arr}    {{\text{\sf Arr}}}
\let\xto \xrightarrow

\let\To \Rightarrow
\let\too \longrightarrow
\let\Too \Longrightarrow
\let\s \boldsymbol
\def\quitar#1{}
\renewcommand{\emph}[1]{{\it#1\/}}
\newcommand{\caja}[2][0ex]{%
\settowidth{\mdpt}{{#2}}%
\settowidth{\mhei}{{#1}}%
\setlength{\mytmplgth}{\mdpt-\mhei}%
#1\hbox to \mytmplgth{}}

\begin{document}
\title{Categorical non abelian cohomology, and the {S}chreier theory
of groupoids}
\author{V.~Blanco, M.~Bullejos, and E.~Faro}

\maketitle
\begin{abstract}
By regarding the classical non abelian cohomology of groups from a 2-dimensional categorical
viewpoint, we are led to a non abelian cohomology of groupoids which continues to satisfy
classification, interpretation and representation theorems generalizing the classical ones. This
categorical approach is based on the fact that if groups are regarded as categories, then, on the
one hand, crossed modules are 2-groupoids and, cocycles are lax 2-functors and the cocycle
conditions are precisely the coherence axioms for lax 2-functors, and, on the other hand group
extensions are fibrations of categories. Furthermore, $n$-simplices in the nerve of a 2-category
are lax 2-functors.
\end{abstract}
\section{Introduction}

In this paper the authors have taken up the task of working out and writing down what has been in
the mind of a few specialists for quite some time. Although this has been a harder job than we
initially imagined, in the end the main point of the paper is making available an \emph{approach}
that may be used by people working in the difficult and artificially disconnected fields of
higher dimensional category theory and higher non abelian cohomology theory.

It is part of the categorical folklore that \emph{``if your groups are not required to be
abelian, you may as well suppose they are groupoids''}. One should think that Schreier's solution
(see \cite{Schreier1926a}) to the classification problem of abelian extensions of groups
---extended by Dedecker \cite{Dedecker1964} to non abelian extensions--- is only part of a larger
story. In the whole story, of course, group extensions would be replaced by some sort of groupoid
extensions, and the non abelian group cohomology containing the ``Schreier invariants'' of non
abelian group extensions, by a more general cohomology of groupoids.

Since a group extension
\begin{equation} \label{gext} 1\to K\xto{j} E\xto{p} G\to 1
\end{equation}
is completely determined by the epimorphism $p$, and since group epimorphisms are the same as
fibrations whose domain and codomain are one-object groupoids, it seems that a categorification
of the theory of Schreier invariants of group extensions would begin with the classification of
fibrations of groupoids.

One such classification is contained in Grothendieck's work since the 2-category, $\fib(\g)$, of
(op)fibrations over any small groupoid $\g$ is 2-equivalent to the 2-category $\LaxFun(\g,\gpd)$
of lax 2-functors from $\g$ to the 2-category \gpd\ of groupoids. In one direction, the
2-equivalence is given by the classical Grothendieck construction \cite{Grothendieck1971} (called
here \emph{twisted product}). An inverse to the Grothendieck construction is obtained (using the
axiom of choice) by associating to any given (op)fibration $\e\to\g$ a lax 2-functor ``fiber''.

The fact that in the study of group extensions such as \eqref{gext} not only the quotient $G$ but
also the kernel $K$ are fixed, indicates that we must classify not all fibrations above a
groupoid $\g$ but only those with fixed fibers. On the other hand, the fact that all group
morphisms, when they are viewed as functors, are bijective on objects make us reduce our
classification of fibrations above a groupoid $\g$ to those which are the identity on objects.
This is equivalent to classifying fibrations above $\g$ whose fibers are groups. Therefore,
``fixing the kernel'' consists in fixing a family of groups $\k$, indexed by the objects of $\g$,
and the problem of classifying all fibrations above $\g$ whose fibers are given by the family of
groups $\k$ translates, by the above 2-equivalence, to classifying not all lax 2-functors from
$\g$ to \gpd\ but only those which have a fixed image. We observe that the family $\k$ defines a
2-groupoid $\Aut(\k)$ (i.e.\ a crossed module), which is a 2-subcategory of \gpd\ and that
fibrations above $\g$ with fibers $\k$ correspond to lax 2-functors which factor through
$\Aut(\k)$. Then, we introduce the category $\Laxact(\g,\k)$ with objects those lax 2-functors
from \g\ to \gpd\ which factor through $\Aut(\k)$ and define a non abelian cohomology of
groupoids (which classifies extensions of groupoids by groups) in terms of connected components
of lax 2-functors
$$
H^2_{\abs{\Aut(\k)}}(\g, \k) = H^1\big(\g,\Aut(\k)\big):=\Laxact[\g,\Aut(\k)]
$$
so that cocycles will arise from a parameterization of lax 2-functors. In this approach
Grothendieck's construction becomes an interpretation theorem of this cohomology in terms of
groupoid extensions (Corollary \ref{inter th}). Furthermore, the \emph{geometric nerve} of
2-categories can be used to give a representation theorem of this cohomology in terms of homotopy
classes of simplicial maps (Theorem \ref{repr th}).

Continuing work within this categorical approach is expected to lead to a non abelian
2-dimensional cohomology of 2-groupoids (very closely related to a non abelian 3-dimensional
cohomology of groups or groupoids) which would classify extensions of 2-groupoids and which would
be represented by homotopy classes of continuous maps from a 2-type to a 3-type.

\section{Weak actions of groupoids and twisted products}
\label{gpdfibrations}

This is the main section of the paper. The relevant facts of 2-dimensional category theory are
presented here with a focus on the concepts (weak actions of groupoids and twisted products) on
which the main results are based. By giving a detailed account of weak actions and fibrations we
are able to make almost all results in the two subsequent sections to appear either evident or as
immediate consequence of the work done here. Although this has made this section to have a
comparatively large size we think the effort has been worthwhile.

\subsection{2-Categories, lax 2-functors, and lax 2-natural
transformations}

A 2-\emph{category} $\ba$ is just a category enriched in the cartesian closed category $\Cat$ of
(small) categories. Thus, $\ba$ consists of 0-cells (or objects), 1-cells (or arrows) and
2-cells, so that for any two objects $X,Y\in \ba$ there is a (small) ``hom'' category $\ba(X,Y)$
whose objects are 1-cells $f:X\to Y$ and whose arrows are 2-cells $\alpha:f\to g$. The
composition in all these categories is globally called the \emph{vertical} composition of $\ba$,
denoted by ``$\circ$'' or, if there is no ambiguity, by simple juxtaposition. The identity map of
$f$ in $\ba(X, Y)$ is denoted $\id{f}$. Connecting the different hom categories there is a
``horizontal composition'' \emph{functor}
\begin{equation} \label{hocomp} \ba(X,Y)\times \ba(Y,Z)\xto{\ *\ }
\ba(X,Z),
\end{equation}
satisfying strict associativity, and in each category $\ba(X,X)$ there is a distinguished object
(i.e.\ arrow of $\ba$), denoted $1_X$, which is a strict right and left identity for the
horizontal composition. Also, for any three objects $X$, $Y$, $Z$, if $f$ is any given 1-cell
$f:X\to Y$, and $\alpha$ is a 2-cell in $\ba(Y,Z)$, we use the customary notation $\alpha*f$ or
simply $\alpha f$ to denote $\alpha*\id{f}$ (and similarly on the other side).

The \emph{underlying category} of $\ba$ will be denoted $\abs{\ba}$. It has as objects those of
$\ba$, as arrows the 1-cells of $\ba$ and its composition is the horizontal composition of
1-cells in $\ba$.

\medskip

An example of a 2-category is provided by $\cat$ itself. Its 2-cells are the natural
transformations, whose vertical composition is essentially given by the composition in the
codomain category.

We will make no notational distinction between $\cat$ regarded as a category and $\cat$ regarded
as a 2-category since the context will always make it clear in what way it is being considered.

\medskip

Another example of 2-category is the full 2-subcategory of $\cat$, denoted $\gpd$, determined by
all small groupoids (categories in which every arrow has an inverse). Also, groups can be
regarded as groupoids having only one object (which will be generically denoted ``$*$'', should
it be necessary to refer to it explicitly), and having the group elements as 1-cells. On the
other hand, group homomorphisms are precisely the functors between groups regarded as groupoids,
so that the category $\Gr$ of groups can be regarded as a full subcategory of $\gpd$. As we do
with $\cat$ and with $\gpd$, we will not make any notational distinction between the category
$\gp$ and the full 2-subcategory of $\gpd$ determined by all groups (of which $\gp$ is the
underlying category). When $\Gr$ is regarded as a 2-subcategory of $\Gpd$, one is implicitly
considering as 2-cells $\alpha:f\to g$ between two group homomorphisms $f,g:G\to H$ those
elements $\alpha\in H$ representing a natural transformation from $f$ to $g$, that is, such that
for every $u\in G$, the following \emph{naturality condition} holds: $g(u)\, \alpha = \alpha
\,f(u).$ Since for any group homomorphism $f:G\to H$ and any element $\alpha\in H$, the map $g$
defined by
\begin{equation} \label{2-cell in Gp} g(u) = \alpha
\,f(u)\,\alpha^{-1}
\end{equation}
is again a group homomorphism, every such pair $(f,\alpha)$ determines a 2-cell in $\gp$ whose
domain is $f$. Its codomain is the homomorphism $g=\alpha f \alpha^{-1}$ defined by \eqref{2-cell
in Gp}. Thus, for any two groups $G$, $H$, the 2-cells in $\gp(G,H)$ are the pairs $(f,\alpha)$
where $f:G\to H$ is a group homomorphism and $\alpha$ is an element of $H$. The vertical
composition is essentially given by the group product in the codomain:
\begin{equation}\label{vert of 2-cell in Gp} \xymatrix@+1pc{G
\ruppertwocell^f{\alpha}\ar[r]|{\;g\;}\rlowertwocell_h{\beta} & H} \qquad (g,\beta)(f,\alpha) =
(f,\beta\alpha),
\end{equation}
while the horizontal composition is a sort of semidirect product,
\begin{equation} \label{horiz of 2-cell in Gp} \xymatrix@+1pc{G
\rtwocell^{f}_{h}{\alpha} & H \rtwocell^{g}_{k}{\;\beta} & K} \qquad (g,\beta)*(f,\alpha) =
\big(gf,\beta g(\alpha)\big).
\end{equation}
When $\gp$ is regarded as a 2-category, any group $G$ determines a full 2-subcategory of $\gp$
which will be denoted $\Endos(G)$. Its underlying category is $\abs{\Endos(G)} = \endos(G)$, the
category (monoid) of endomorphisms of $G$.

\bigskip

A 2-\emph{groupoid} $\ba$ is a 2-category whose underlying category as well as each of its hom
categories are groupoids. This implies that 2-cells are invertible, not only for the vertical
composition (for which the inverses are denoted $\alpha^{-1}$), but also for the horizontal
composition, the horizontal inverse of a 2-cell being given by the formulas
\begin{equation} \label{two defs of inv} \inv(\id f) = \id{f^{-1}}\,,
\qquad \inv(\alpha) = f^{-1} * \alpha^{-1} * g^{-1} = g^{-1} * \alpha^{-1} * f^{-1},
\end{equation}
(where $\alpha:f\to g$ is a 2-cell).

In the category $\gp$ of groups not every 1-cell (group homomorphism) has an inverse. However,
all hom categories $\gp(G,H)$ are groupoids since every 2-cell $(f,\alpha)$ in $\gp$ has an
inverse given by
$$
(f,\alpha)^{-1} = \big(\alpha f \alpha^{-1},\alpha^{-1}\big).
$$
Thus, if we restrict the 1-cells in $\gp$ to the group isomorphisms, we obtain a sort of (large)
2-groupoid which will be denoted $\iso(\gp)$.

In the same way that a group $G$ determines a full 2-subcategory of $\gp$, it also determines a
full 2-subcategory of $\iso(\gp)$, denoted $\Aut(G)$, whose underlying category is $\abs{\Aut(G)}
= \aut(G)$, the group of automorphisms of $G$. Note that $\Aut(G)$ is automatically a 2-groupoid.
We will find it useful to generalize the above notation to arbitrary families of groups, so that
if $\k = \{K_x\}_{x\in X}$ is a family of groups indexed by a set $X$, the full 2-subcategory of
$\iso(\gp)$ determined by $\k$ will be denoted $\Aut(\k)$. As in the previous case, this is
actually a 2-groupoid.

\bigskip

A 2-\emph{functor} $F:\ba\to \bb$ is an enriched functor in $\Cat$, so it takes objects, 1-cells
and 2-cells in $\ba$ to objects, 1-cells and 2-cells in $\bb$ respectively, in such a way that
all the 2-category structure of $\ba$ is strictly preserved. In particular, any 2-functor
preserves inverses of 1-cells and of 2-cells.

Small 2-categories and 2-functors form a category that we denote $\doscat$ (this is actually a
2-category and even a 3-category, but we will not use in this paper those higher dimensional
structures). The full subcategory of $\doscat$ determined by the 2-groupoids is denoted by
$\dosgpd$.

\begin{example}[2-functors $G\to\Aut(K)$] \label{acciones} Given any
two groups $G$ and $K$ we can consider 2-functors $G\to\Aut(K)$. Since $G$ is 2-discrete as a
2-groupoid, to give one such 2-functor is the same as giving a group homomorphism $G\to\aut(K)$,
that is, a group action of $G$ on $K$.
\end{example}

For any small category $\c$, the category of all actions of $\c$ on groups is defined as the
functor category $\gp^\c$. If $\c$ has only one object an action of $\c$ on groups (regarded as a
functor $F:\c\to\gp$) determines the specific group on which $\c$ acts as the image $F(*)$ of the
only object of $\c$. This allows us to parameterize the actions $\c$ by the groups acted upon.
This is one of the ways in which one can arrive at the concept of an action of $\c$ on a specific
group $H$. If we drop all restriction on the number of objects of $\c$, an action of $\c$ on
groups determines, in general, not a group, but \emph{a family of groups indexed by the objects
of $\c$}. This allows us to parameterize the actions of $\c$ on groups by the $\obj{\c}$-indexed
families of groups on which $\c$ may act. (Throughout the paper we will use the notation
$\obj{\c}$ to denote the set of objects of a small category $\c$.)

\begin{definition} If \g\ is a groupoid and \k\ is a family of groups
indexed by the objects of \g, a groupoid action of \g\ on \k\ is a 2-functor $\g\to\Aut(\k)$ (or,
equivalently, a functor $\g\to\abs{\Aut(\k)}$) whose objects function is the indexing of $\k$,
that is, such that for every object $A\in\g$, $K_A = F(A)$.
\end{definition}

\bigskip

We come now to the maps between 2-categories that are crucial in this paper. These are the
\emph{lax} 2-functors, whose definition we recall next:

\begin{definition}\label{def lax 2-funct}

Given two 2-categories $\ba$, $\bb$, a (normal\,\footnote{Note that in this paper all lax
2-functors will be normalized. Thus, the expression ``lax 2-functor'' will mean here what is
usually referred to in the literature as ``normal lax 2-functor'', that is, for us a lax
2-functor is strictly identity preserving. The structural natural transformation for identities
is omitted from the definition since it is itself an identity. In exchange, we need to include an
additional requirement of ``coherence with the identity law'' for the structural natural
transformation for composition, $\sigma$.}) lax 2-functor from $\ba$ to $\bb$ is a pair
$(F,\sigma)$ where

\begin{myitemize}

\item[(a)] $F$ is a correspondence which takes each object $X$ of $\ba$ to
an object $F(X)$ of $\bb$, and any two objects $X$, $Y$ in $\ba$, to a functor $F_{XY} : \ba(X,Y)
\to \bb\big(F(X),F(Y)\big)$.
\item[(b)] $\sigma$ is a correspondence which assigns to any three objects
$X$, $Y$, $Z$ in $\ba$ a natural transformation $\sigma^{XYZ}$ between functors from
$\ba(X,Y)\times\ba(Y,Z)$ to $\bb\big(F(X),F(Z)\big)$, whose component on an ``object'' $X\xto f Y
\xto g Z$ is a 2-cell in $\bb$ denoted $\sigma^{XYZ}_{gf}: F_{XZ} (gf) \to F_{YZ} (g) F_{XY}
(f)$.
\end{myitemize}

We will generally omit the subscripts in $F_{XY}$ and the superscripts in $\sigma^{XYZ}$, the
objects being implicit, so that the above component of $\sigma$ will be written
\begin{equation}\label{es dos celda} \sigma_{gf}: F(gf) \to F(g) F(f).
\end{equation}
The naturality of $\sigma$ gives, for each pair of horizontally composable 2-cells
$$
(\alpha, \beta):(f,g)\to (h,k),\quad \text{or} \quad \xymatrix{\cdot \rtwocell^f_{h}{\alpha}
&\cdot \rtwocell^g_{k}{\beta} &\cdot},
$$
the following equation:
\begin{equation}\label{nat condit} \sigma_{kh} \,F(\beta*\alpha) =
\big(F(\beta) * F(\alpha)\big) \, \sigma_{gf}.
\end{equation}
$$
\vcenter{ \xymatrix{(f,g) \ar[d]_{(\alpha,\beta)} \\ (h,k)} } \quad \vcenter{ \xymatrix{ F(gf)
\ar[d]_{F(\beta*\alpha)}
\ar[r]^-{\sigma_{gf}} &F(g)F(f) \ar[d]^{F(\beta) * F(\alpha)} \\
F(kh)\ar[r]_-{\sigma_{kh}} &F(k)F(h)} }.
$$
These data are required to satisfy the following axioms:

\begin{myenumerate}

\item[LF1.] (Normalization) $F(1_X) = 1_{F(X)}$.

\item[LF2.] (Coherence with associative law) For any 3-chain in $\ba$,
$X\xto{f}Y\xto{g}Z\xto{h}T$,
\begin{equation}\label{ax coher assoc}
\big(\sigma_{hg}*F(f)\big)\,\sigma_{(hg)f}\, = \big(F(h)*\sigma_{gf}\big)\, \sigma_{h(gf)}\, ,
\end{equation}
$$
\vcenter{ \xymatrix @C=3pc{ F(hgf)\ar[d]_{\sigma_{h(gf)}}
\ar[r]^-{\sigma_{(hg)f}} &F(hg)F(f) \ar[d]^{\sigma_{hg}*F(f)} \\
F(h)F(gf) \ar[r]_-{F(h)*\sigma_{gf}} & F(h)F(g)F(f).} }
$$
\item[LF3.] (Coherence with identity law) For any 1-cell $f:X\to Y$,
$$
\sigma_{f1_X} = \id{F(f)} = \sigma_{1_Yf}.
$$
\end{myenumerate}

\end{definition}

Small 2-categories and lax functors form a category that we denote by $\ldoscat$, analogously
$\ldosgpd$ will be the full subcategory of $\ldoscat$ whose objects are the 2-groupoids.

\begin{proposition}\label{weak id coherence} In the presence of
inverses of 2-cells in the codomain 2-category, axiom LF3 in Definition \ref{def lax 2-funct} can
be substituted by the following weaker version:
\begin{myenumerate}
\item[LF3'.] (Weak coherence with identity law) For any object $X$,
$$
\sigma_{1_X1_X} = \id{F(1_X)}( = \id{1_{F(X)}}).
$$
\end{myenumerate}
\end{proposition}

\begin{proof} Putting $g=h=1_Y$ in \eqref{ax coher assoc} and using
axiom LF1 and the inverse 2-cell of $\sigma_{1_Y f}$, one gets $\sigma_{1_Yf} = \sigma_{1_Y
1_Y}*F(f) = \sigma_{1_Y 1_Y}*\id{F(f)}$. Therefore axiom LF3' implies $\sigma_{1_Yf} =
\id{F(f)}$. Similarly, putting $f=g=1_Z$ in \eqref{ax coher assoc} and using axiom LF1 and the
inverse 2-cell of $\sigma_{h 1_Z}$, one gets $\sigma_{h 1_Z}*F(1_Z) = F(h)*\sigma_{1_Z 1_Z} =
\id{F(h)}*\sigma_{1_Z 1_Z}$. Therefore axiom LF3' implies $\sigma_{h 1_Z} = \id{F(h)}$.
\end{proof}

Since we have enriched the group $\aut(K)$ of automorphism of a group $K$ to the 2-category
$\Aut(K)$, in such a way that a group action of a group $G$ on $K$ is just a 2-functor
$G\to\Aut(K)$, we can relax the definition of action by considering lax 2-functor $G\to \Aut(K)$
which are just weak actions of $G$ on $K$. Moreover there is no difficulty to passing from group
to groupoids.

\begin{definition}[Weak groupoid actions] \label{def lax act} If $\g$
is a groupoid and $\k = \{K_A\}_{A\in\g}$ is a family of groups indexed by the objects of $\g$, a
\emph{weak action} of $\g$ on $\k$ is a lax 2-functor $F:\g\to \Aut(\k)$ whose objects function
is the indexing of $\k$, that is, such that for every object $A\in\g$, $K_A = F(A)$.
\end{definition}

If $G$ is a group, the condition on the lax 2-functor $G\to\Aut(K)$ to be a lax action is always
satisfied and therefore a lax group action of $G$ on $K$ is just a lax 2-functor $G\to\Aut(K)$.

The following proposition spells out the data and the axioms characterizing weak actions of
groupoids.

\begin{proposition} \label{cocycle} Given a groupoid \g\ and a family
of groups \k\ indexed by the objects of $\g$, a weak action of \g\ on $\k$ is a pair $\langle F,
\sigma\rangle$ where $F$ is a correspondence assigning to each arrow $u:A\to B$ in $\g$ a group
isomorphism $F(u) = \lact u{\,(-)}:K_A \to K_B$, and $\sigma$ is a correspondence assigning to
each composable pair $A\xto u B \xto v C$ in $\g$ a group element $\sigma_{vu}\in K_C$ such that,
for each 3-chain $A\xto u B \xto v C \xto w D$ \ in $\g$, and any element $x\in K_A$, the
following equations are satisfied,
\begin{myenumerate}
\item $\lact{1}{\,x} = x$,
\item $\sigma_{11} = 1$,
\item $\sigma_{vu}\lact{vu}{\,x} = \lact v{\big(\lact
u{\,x}\big)}\sigma_{vu}$,
\item $\sigma_{wv}\sigma_{(wv)u} = \lact w{\,\sigma_{vu}}
\sigma_{w(vu)}$.
\end{myenumerate}
\end{proposition}
\begin{proof} Condition 3 can be rewritten $F(v)\circ F(u) =
\sigma_{vu} F(vu) \sigma_{vu}^{-1}$ and it assures, according to equation \eqref{2-cell in Gp}
(page \pageref{2-cell in Gp}), that $\sigma_{vu}$ is a 2-cell from $F(vu)$ to $F(v)\circ F(u)$.
Since $\g$ is 2-discrete, the naturality condition \eqref{nat condit} is vacuous and therefore
$\sigma$ is a natural transformation. Condition 1 is the normalization of $F$, condition 2, the
coherence with the identity law, and condition 4, the coherence with the associative law.
\end{proof}

If $G$ is a group, a weak action of $G$ on a group $H$ consists of a pair of maps:
$$
F:G\to\aut(K), \qquad \sigma:G\times G\to K
$$
satisfying the conditions in Proposition \ref{cocycle} (in this context $\sigma(u,v)$ will denote
the element $\sigma_{vu}$). This is just a Dedecker's 2-cocycle of $G$ with coefficients on the
crossed module $K\to\aut(K)$ associated to the 2-groupoid $\Aut(K)$.

The natural morphisms between lax 2-functors are the lax 2-natural transformations, whose
definition we recall next.

\begin{definition} Given lax 2-functors $F_1,F_2:\ba\to \bb$, a lax
2-natural transformation from $F_1$ to $F_2$ is a pair $\lalpha = (\alpha, \tau^\alpha)$ where
\begin{myitemize}
\item[(a)] $\alpha$ is a correspondence which assigns to each object $A\in
\ba$ an arrow $\alpha_A:F_1(A)\to F_2(A)$ in $\bb$, and
\item[(b)] $\tau^\alpha$ is a correspondence which assigns to each two
objects $A$, $B$ in $\ba$ a natural transformation between functors from $\ba(A,B)$ to
$\bb\big(F_1(A), F_2(B)\big)$, whose component at a morphism $f:A\to B$ in $\ba$ is a 2-cell,
$\tau^\alpha_f$, in $\bb$ as in the diagram
$$
\vcenter{ \xymatrix @C=.7pc @R=2pc { F_1(A) \ar[d]_{F_1(f)} \ar[rr]^-{\alpha_{A}}&
\ar@{}[d]|-(.39){\overset{\scriptstyle\tau^\alpha_f}{\textstyle\Too} } & F_2(A) \ar[d]^{F_2(f)}
\\ F_1(B) \ar[rr]_-{\alpha_{B}} & & F_2(B) }}\qquad \tau^\alpha_f : \alpha_{B}\,F_1(f) \to
F_2(f)\,\alpha_{A},
$$
the naturality of $\tau^\alpha$ gives, for each arrow $\eta:f\to g$ in $\ba(A,B)$ the following
equation:
$$
\big(F_2(\eta)*\alpha_{A}\big)\,\tau^\alpha_f = \tau^\alpha_g\,\big(\alpha_{B}*F_1(\eta)\big).
$$
\end{myitemize}

These data are to satisfy the axioms:

\begin{myenumerate}
\item (normalization) $\tau^\alpha_{1_A}=\id{{\alpha_A}}$,
\item (coherence) for each composable pair of arrows
${A}\xto{f}{B}\xto{g}{C}$,
$$
\big(F_2(g)*\tau^\alpha_f\big) \big(\tau^\alpha_g*F_1(f)\big)
\big(\alpha_{C}*{\sigma_1}_{gf}\big) = \big({\sigma_2}_{gf}*\alpha_{A}\big)\tau^\alpha_{(gf)}.
$$
\end{myenumerate}
\end{definition}

An example of lax 2-natural transformation is given by the \emph{identity lax 2-natural
transformation} of a lax 2-functor $F:\ba\to\bb$, which is defined by $\alpha_A = 1_{F(A)}$ and
$\tau_f = \id{F(f)}$.

Given two lax 2-natural transformations $F_1\xto{\lalpha}F_2 \xto{\lbeta} F_3$ between lax
2-functors $\ba\to\bb$, a composite lax 2-natural transformation, $\lbeta\lalpha: F_1\to F_3$,
can be defined with structure given by
$$
\tau^{\beta\alpha}_f = (\tau^\beta_f * \alpha_A) (\beta_B * \tau^\alpha_f).
$$
This composition is associative and it has the identity lax 2-natural transformations as
two-sided identities so that for any two 2-categories $\ba$, $\bb$ there is a category, denoted
$\LaxFun(\ba,\bb)$, whose objects are the lax 2-functors $\ba\to\bb$ and whose arrows are the lax
2-natural transformations between them. By taking lax 2-natural transformations as 2-cells, the
category $\ldoscat$ of small 2-categories and lax 2-functors becomes a 2-category.

\begin{example}[Maps in $\LaxFun\big(\g,\Aut(\k)\big)$] \label{lax
2-nat script G to a fam of gps} Let $\g$ be a groupoid, $\k = \{K_A\}_{A\in\g}$ a family of
groups indexed by the objects of $\g$, and let $(F_1,\sigma^1),(F_2,\sigma^2):\g\to\Gr$ be lax
2-functors from $\g$ to $\Aut(\k)$. A lax 2-natural transformation $\lalpha:F_1\to F_2$ is a pair
$\lalpha = (\alpha, \tau)$ where for every object $A\in\g$, $\alpha_A: K_A \to K_A$ is a group
isomorphism, and for every pair of objects $A,B\in \g$, $\tau_{AB}: \g(A,B) \to K_B$ is a map of
sets satisfying for all $ k\in K_x$, and $u\in \g(A,B)$,

\begin{myenumerate}

\item \caja[(naturality)]{(normalization)}
$F_2(u)\big(\alpha_A(k)\big) \cdot \tau_{AB}(u) = \tau_{AB}(u)\cdot\alpha_B \big(F_1(u) (k)\big)
$,

\item (normalization) $\tau_{AA}(1) = 1$, and

\item \caja[(coherence)]{(normalization)} $
F_2(v)\big(\tau_{AB}(u)\big) \cdot \tau_{BC}(v) \cdot \alpha_C(\sigma_{vu}^1) = \sigma_{vu}^2
\cdot \tau_{AC}(vu)$
\end{myenumerate}
\end{example}

We next define the category of weak actions of a groupoid $\g$ on a family $\k$ of groups indexed
by the objects of $\g$.

\begin{definition} \label{def lax act} Let $\g$ be a groupoid, and $\k
= \{K_A\}_{A\in\g}$ a $\obj{\g}$-indexed family of groups. If $F_1$, $F_2$ are two weak actions
of $\g$ on $\k$, a map of weak actions $F_1\to F_2$ in the category $\Laxact(\g,\k)$ is a lax
2-natural transformations $(\alpha,\tau):F_1 \to F_2$ such that for every object $A\in\g$,
$\alpha_A = 1_{K_A}$.
\end{definition}

According to this definition, there is an inclusion of categories
$$
\Laxact(\g,\k)\subseteq\LaxFun\big(\g,\Aut(\k)\big)
$$
which is an identity on objects. In the next example we see how the axioms of lax 2-natural
transformations get simplified in the case of weak actions of a group.

\begin{example}[Morphisms of weak actions of a group] \label{arrows
between weak actions} Let $G$ and $K$ be groups and let
$(F_1,\sigma_1),(F_2,\sigma_2):G\to\Aut(K)$ be weak actions of $G$ on $K$. To give a morphism
$F_1\to F_2$ in $\Laxact(G,K)$ is to give a map of sets $\tau: G \to K$ satisfying for all $ k\in
K$, and $u,v\in G$,

\begin{myenumerate}
\item \caja[(naturality)]{(normalization)} $\tau(u)\cdot F_1(u) (k) =
F_2(u)(k) \cdot \tau(u)$,
\item (normalization) $\tau(1) = 1$,
\item \caja[(coherence)]{(normalization)} $F_2(v)\big(\tau(u)\big)
\cdot \tau(v) \cdot \sigma_1(u,v) = \sigma_2(u,v) \cdot \tau(vu)$
\end{myenumerate}

Note that by the naturality condition, $F_2$ is completely determined by $\tau$ and $F_1$ and
that the coherence condition determines $\sigma_2$:
\begin{align}\label{solve for F} F_2(u)(k) &= \tau(u)\cdot F_1(u) (k)
\cdot \tau(u)^{-1} \\ \label{solve for sigma} \sigma_2(u,v) &= \tau(v)\cdot F_1(v)
\big(\tau(u)\big) \cdot \sigma_1(u,v) \cdot \tau(vu)^{-1}.
\end{align}

If we write $\mathit{C}^1(G,K)$ for the set of maps $\tau:G\rightarrow K$ such that $\tau(1)=1$
and $\mathit{C}^2_{\aut(K)}(G)$ for the set of arrows in $\Laxact(G,K)$ we can define an
operation
$$
\nabla: \mathit{C}^1(G,K)\times\mathit{C}^2_{\aut(K)}(G) \too \mathit{C}^2_{\aut(K)}(G)
$$
as $\nabla(\tau,F_1)=\tau\nabla F_1=(F_2,\sigma_2)$ where $F_2$ y $\sigma_2$ are defined by
equations \eqref{solve for F} and \eqref{solve for sigma}. The map $\nabla$ is one of the three
operations defined by Dedecker in \cite{Dedecker1964} where it is used to define cohomologous
cocycles.
\end{example}

\begin{example}[Morphisms of weak actions of a groupoid] \label{w.act
map script G to a fam of gps} Let $\g$ be a groupoid, let $\k = \{K_A\}_{A\in\g}$ be a family of
groups indexed by the objects of $\g$, and let $(F_1,\sigma^1),(F_2,\sigma^2):\g\to\Gr$ be weak
actions of $\g$ on $\k$. To give a morphism $F_1\to F_2$ in $\Laxact(\g,\k)$ is to give for every
pair of objects $A,B\in \g$, a map of sets $\tau_{AB}: \g(A,B) \to K_B$ satisfying, for all $
k\in K_A$, and $A\xto u B \xto v C$ in $\g$,

\begin{myenumerate}
\item \caja[(naturality)]{(normalization)} $\tau_{AB}(u)\cdot F_1(u)
(k) = F_2(u)(k) \cdot \tau_{AB}(u) $,
\item (normalization) $\tau_{AA}(1) = 1$,
\item \caja[(coherence)]{(normalization)} $
F_2(v)\big(\tau_{AB}(u)\big) \cdot \tau_{BC}(v) \cdot \sigma_{vu}^1 = \sigma_{vu}^2 \cdot
\tau_{AC}(vu)$
\end{myenumerate}

As in the group case, the codomain of a morphism of weak actions of $\g$ on $\k$ is completely
determined by its source and the family of maps $\tau_{AB}$. Then, if we write
$\mathit{C}^1(\g,\k)$ for the set of families of maps
$\boldsymbol{\tau}=\{\tau_{AB}\}_{A,B\in\g}$ such that $\tau_{AA}(1)=1$ for all object $A\in \g$,
and $\mathit{C}^2_{\Aut(\k)}(\g)$ for the set of arrow in $\Laxact(\g,\k)$, we can define an
operation
$$
\nabla: \mathit{C}^1(\g,\k)\times\mathit{C}^2_{\Aut(\k)}(\g) \too \mathit{C}^2_{\Aut(\k)}(\g)
$$
as $\nabla(\boldsymbol{\tau},F_1)=\boldsymbol{\tau} \nabla F_1=(F_2,\sigma^2)$ where $F_2$ y
$\sigma^2$ are obtained by solving for them in the naturality and coherence conditions above.
\end{example}

\subsection{Fibrations and the lax 2-functor ``fiber'' associated to a
fibration}

Given a functor, $F:\a\to\g$, the fiber of $F$ over an object $B\in \g $ is the subcategory
$F^*(B)$ of $\a$ defined by the following pullback in $\cat$:
\begin{equation}\label{pullback fiber} \vcenter{ \xymatrix@=5mm{
F^*(B) \ar[d] \ar[r] \ar@{}[dr]|-(.45){\scriptscriptstyle\text{p.b.}} & \s1 \ar[d]^{B} \\ \a
\ar[r]_{F} & \g } }
\end{equation}
Thus, $F^*(B)$ has:
\begin{myitemize}
\item objects: those objects $A\in \a$ such that $F(A)=B$, and
\item arrows: those arrows $f\in\a$ such that $F(f)=\id{B}$.
\end{myitemize}

A fibration is a functor for which the process of ``taking fibers'' is as functorial as it can
be.

\begin{definition} A functor between groupoids $F:\a\to \g $ is called
a (Grothendieck op-) fibration if for any object $A\in \a$ and any arrow $f:F(A)\to B$ in $\g $
there exists an arrow $\til{f}:A\to \til{A}$ in $\a$ such that $F\big(\til{f}\,\big)=f$.
\end{definition}

We will write $\fib(\g )$ for the full subcategory of the slice category $\Gpd/\g $ determined by
those objects which are fibrations.

\begin{example}[The lax 2-functor ``fiber'' of a fibration of
groupoids] \label{fiber} Given a fibration $F:\a\to\g $ of groupoids, let us suppose we have
chosen for each pair $(f,A)$ \big(where $f:X\to Y$ is an arrow of $\g $ and $A$ is an object of
$F^*(X)$\big), an arrow $\charr{f}{A}$ with domain $A$ and such that $F(\charr{f}{A})=f$ (if $f$
is an identity, the choice is $(\til{1_{X}})_A=1_{A}$). Let's denote by $\lact{f}A$ the codomain
of $\charr{f}{A}$ so that we have $\charr{f}{A}:A\to\lact{f}{A}$, and
$$
\xymatrix@R=4ex@C=2ex{A \ar@{.>}[d]_-{\charr{f}{A}} &
\ar@{}[d]|-(.4){\textstyle\stackrel{F}{\mapsto}} & X\ar[d]^f \\
\lact{f} A & & Y}
$$
Based on the above choices one can give a lax 2-functor ``fiber of $F$'' (which is not unique,
since it depends on a choice),
$$
(F^*,\sigma):\g \to \Gpd,
$$
defined in the following way:
\begin{myitemize}
\item[$F^*$:] On objects, via the pullback \eqref{pullback fiber}; on
the arrows $f:X\to Y$ in $\g $ as the functor $F^*(f):F^*(X) \to F^*(Y)$ defined by (where
$h:A\to B$ in $F^*(X)$)
\begin{equation}\label{def F star of f} F^*(f)(A) = \chobj{f}A\,,
\qquad F^*(f)(h) = \charr{f}{B} h {\charr{f}{A}}^{-1}.
\end{equation}
\item[$\sigma$:] Its component on a composable pair $X \xto f Y \xto g
Z$ in $\g $ is the natural transformation $\sigma_{gf}: F^*(gf) \to F^*(g)F^*(f)$ whose component
at an object $A$ of $F^*(X)$, is the arrow $\sigma_{gf}^A: \chobj{gf}{A} \to
\lact{g}{\big(\lact{f}{A}\big)}$, defined by
$$
\sigma_{gf}^A = \charr{g}{(\lact[-.3ex]{f}{A})}\; \charr{f}{A}\;{\charr{(gf)}A}^{-1}.
$$
\end{myitemize}
\end{example}

It is a simple exercise to prove that equations \eqref{def F star of f} define a functor
$F^*(f)$, that $\sigma_{gf}$ is indeed a natural transformation, and that $\sigma$ satisfies the
conditions for the structure map of a lax 2-functor (note that in this context one gets the
equation
\begin{equation}\label{coh in gpd}
\sigma_{hg}^{\lact{f}{A}}\;\sigma_{(hg)f}^A = F^*(h)\big(\sigma_{gf}^A\big)\;\sigma_{h(gf)}^A\,,
\end{equation}
for the associativity coherence of $\sigma$), so that the above data $(F^*,\sigma)$ indeed
determines a lax 2-functor.

\begin{observation}\label{isos in fibrations} It is noteworthy the
fact that for each arrow $f:X\to Y$ in $\g$, the functor $F^*(f)$ is bijective at the level of
arrows. That is, for any two objects $A, B\in F^*(X)$ the map
$$
F^*(f):F^*(X)(A,B) \too F^*(Y)(\chobj{f}A,\chobj{f}B)
$$
defined by \eqref{def F star of f} is a bijection (the inverse map sends $k$ to
${\charr{f}{B}}^{-1} k \charr{f}{A}$). In particular, for any object $A\in\a$ and any arrow
$f:X\to Y$ in $\g$, $F^*(f)$ determines a group isomorphism between the kernels of the group
homomorphisms $F_{AA}:\a(A,A)\to\g\big(X,X\big)$ and $F_{\chobj{f}A \chobj{f}A}: \a(\chobj{f}A,
\chobj{f}A) \to \g\big(Y,Y\big)$.
\end{observation}

\begin{example}[Case of a fibration bijective on objects] \label{fiber
of a fibration b.o.} If a fibration of groupoids $F:\e\to\g$ is bijective on objects, every fiber
has only one object and it is therefore a group, so that any fiber functor $F^*$ of $F$ goes to
$\gp$. Furthermore, by Observation~\ref{isos in fibrations}, every arrow of $\g$ is sent by $F^*$
to a group isomorphism. The data that needs to be given together with $F$ in order to determine a
lax functor ``fiber'' $F^*:\g\to\Gpd$ consists in choosing, for each arrow $f\in \g$, an arrow
$\til{f}\in \e$ such that $F\big(\til{f}\,\big) = f$, and $\til 1=1$, that is, we need to specify
a set-theoretic section of the arrows function of $F$ such that it takes the identities to
identities. Once we have made this choice, the lax 2-functor $F^\ast:\g\to \iso(\gp)$ takes an
object $A$ of $\g$ to the kernel $K_A$ of the group homomorphism
$F_{AA}:\e(A,A)\to\g(A,A)$,\footnote{It is harmless to assume that the (bijective) object
function of $F$ is an identity.} each arrow $f:A\to B$ in $\g$ to the group \emph{isomorphism}
$$
F^\ast(f)=\lact{f}{\,(-)}:K_A\to K_B\, ; \quad \lact{f} {\,k}=\til{f}\,k\,\til{f}^{-1},
$$
and the structure map of $F^*$ is given by
$$
\sigma_{gf} = \til{g}\,\til{f}\;\big(\til{g f}\big)^{-1} .
$$
\end{example}

Is $F$ is a fibration of groupoids which is bijective on objects, any lax 2-functor fiber $F^*$
determines a \emph{weak action} of the base groupoid on the family of groups determined by the
fibers. In the particular case of groups this reduces to the fact that any section of an
epimorphism of groups determines a weak action of the codomain on the kernel.

\begin{example}[The lax 2-functor fiber of a fibration of groups]
\label{fiber of a fibration of groups}

A fibration between groups is just a surjective homomorphism $p:E\twoheadrightarrow G$. The fiber
of $p$ over the only object of $G$ is obviously the kernel of $p$. A lax 2-functor ``fiber''
$p^\ast:G\to \Aut(K)$ is determined by choosing, for each element $u\neq 1$ in $G$, an element
$\til{u}\in E$ such that $p(\til{u})= u$, or by specifying a set-theoretic section of $p$ taking
the identity of $G$ to the identity of $E$. Having made this choice, $p^\ast$ is defined as
taking the only object of $G$ to $K=\ker(p)$, each element $u\in G$ to the following automorphism
of $K$:
$$
p^\ast(u)=\lact{u}{\,(-)}:K\to K\, ; \quad \lact{u} {\,x}=\til{u}\,x\,\til{u}^{-1},
$$
and the structure map of $p^*$ is given by
$$
\sigma:G\times G\to K\, ; \qquad \sigma(u,v)= \sigma_{vu}= \til{v}\,\til{u}\;\til{vu}^{-1}.
$$
\end{example}

\subsection{Twisted Products}

We have seen in Example~\ref{fiber} how fibrations of groupoids give rise to lax 2-functors to
the 2-category $\Gpd$. The inverse process produces a fibration of groupoids from a lax 2-functor
\begin{equation} \label{lax func G to Gpd} F:\g \to \gpd.
\end{equation}
Using the terminology of Schreier \cite{Schreier1926a}, one can call the resulting fibration the
\emph{twisted product} of $\g$ and $F$, since it reduces to the case dealt with by Schreier in
the case that $\g$ is a group, $H=F(*)$ is also a group. The domain of the obtained fibration is
(in the case of groups) the group twisted product of $\g$ and $H$ relative to the weak action
induced by $F$. Although we are using Schreier terminology we will use Grothendieck notation, so
that we write $\semi F\g$ for the domain of the obtained fibration. Let's review this
construction in the particular case at hand.

Given the lax 2-functor \eqref{lax func G to Gpd}, if $f:X\to Y$ is an arrow in $\g$, the action
of the functor $F(f):F(X)\to F(Y)$ will be denoted by a left action notation, $\lact f{(-)}$.
Accordingly, for any arrow $\lambda:A\to B$ in $F(X)$, we write $F(f)(\lambda)=
\big(\lact{f}\lambda:\lact{f} A \to \lact{f}B\big)$.

\begin{definition}[Twisted product] Given a groupoid $\g$ and a lax
2-functor $(F,\sigma):\g \to \gpd$, the \emph{twisted product} groupoid $\semi F\g$ has as
objects the pairs $(X,A)$ where $X$ is an object in $\g$ and $A$ is an object in $F(X)$. The
arrows of $\semi F\g$ are also pairs $(f,\lambda):(X,A)\to (Y,B)$ where $f:X\to Y$ is an arrow in
$\g$ and $\lambda: \lact f{\kern-.06em A} = F(f)(A)\to B$ is an arrow in $F(Y)$. The composition
of arrows
$$
\xymatrix{(X,A)\ar[r]^{(f,\lambda)} & (T,B)\ar[r]^{(g,\mu)} & (Z,C)}
$$
is defined as
$$
(g,\mu)\; (f,\lambda) = \big(gf,\mu \; \lact{g\,}{\lambda} \; \sigma_{gf}^A\big)
$$
$$
\def\labelstyle{\textstyle} \xymatrix@C=1.8pc {\lact{(gf)}{\kern-.1em
A} \ar@<-.5ex>[r]^-{\sigma_{gf}^A} &\lact{g}{\big(\lact{f}{\kern-.1em A}\big)}
\ar@<-.5ex>[r]^-{\lact{g\,}{\lambda}\,} &\lact[1pt]{g}{B} \ar@<-.5ex>[r]^-{\mu} & C},
$$
while the identity map of $(X,A)$ is $(1_X,1_A)$. This gives a category $\semi F \g$ which is a
groupoid since every arrow $(f,\lambda)$ has an inverse, which is defined by
$$
(f,\lambda)^{-1} = \Big(f^{-1}, \big(\lact[1pt]{f^{-1}}{\kern-.1em\lambda} \;
\sigma_{f^{-1}f}^A\big)^{-1}\Big).
$$
(That this is a \emph{right} inverse to $(f,\lambda)$ is the non trivial part of the proof. It
requires the identity $\lact{f}{\big(\sigma_{f^{-1}f}^A\big)} =
\sigma_{ff^{-1}}^{\lact{f}{\kern-.1em A}}$ which is a consequence of \eqref{coh in gpd} with
$h=f$ and $g=f^{-1}$.)

The obvious projection $\semi F \g\to \g$ is a fibration of groupoids, with fiber over any object
$X\in\g$ the groupoid $F(X)$ (after the opportune identification of the objects $A\in F(X)$ with
the objects $(X,A)\in \semi F \g$ and similarly for the arrows).
\end{definition} The twisted product construction just defined is
moreover functorial, that is, one can define, for each lax 2-natural transformation
$\lalpha:F_1\to F_2$ in $\LaxFun(\g,\Gpd)$ a \emph{functor of fibrations} above $\g$,
$$
\textstyle \semi{\lalpha}\g: \semi{F_1}\g \too \semi{F_2}\g
$$
so that given $F_1 \xto{\lalpha} F_2 \xto{\lbeta} F_3$, $\semi{(\lbeta\lalpha)}\g =
\big(\semi{\lbeta}\g\big) \big(\semi{\lalpha} \g\big)$.

So, we get a functor
\begin{equation} \label{equiv} \semi {} \g :\LaxFun(\g,\Gpd)\too
\fib(\g)
\end{equation}

Let $F:\e\to \g$ be a fibration of groupoids. Then, applying the twisted product construction to
any fiber lax functor $F^\ast:\g\to \Gpd$ obtained from $F$ produces a new fibration $\semi
{F^\ast} \g \to \g$ such that there is a functor $\Gamma :\e\to \semi {F^*} \g$ which makes
commutative the triangle
$$
\xymatrix@C=0mm@R=5mm{ **[l]{\phantom{\semi{F^\ast} \g}\e} \ar[rr]^-{\Gamma } \ar[dr]_{F} & &
**[r]{\semi{F^\ast} \g\phantom{\e}} \ar[dl] \\ & \g }
$$
that is, $\Gamma$ is an arrow in $\fib(\g)$.

Furthermore $\Gamma$ is an isomorphism of categories, so that any fibration of groupoids can be
recovered (up to isomorphism) as the canonical projection from a ``\emph{twisted product}". We
have the following important well known result:

\begin{proposition} \label{clasific fibr de gpds} \label{lax act are
ext} The functor $\semi{}\g :\LaxFun(\g,\Gpd)\to \fib(\g)$ is an equivalence of categories.
\end{proposition}

\begin{example}[Twisted product of weak actions] \label{ex lax act and
fib b.o.} For a lax 2-functor $F: \g \to \Gp$, the twisted product $\semi F \g$ is a groupoid
with the ``same objects'' as $\g$ and the twisted product fibration $\semi{F} \g \to \g$ is the
identity on objects. Conversely, if $F:\e\to\g$ is a fibration of groupoids which is bijective on
objects, then any fiber lax 2-functor of $F$ goes to groups. Therefore, the equivalence of
categories
$$
\semi {} \g :\LaxFun(\g ,\Gpd)\too \fib(\g )
$$
restricts to an equivalence of categories
$$
\semi {} \g : \LaxFun(\g ,\Gr) \too \bofib(\g ),
$$
where $\bofib(\g )$ is the full subcategory of $\fib(\g )$ determined by those fibrations which
are bijective on objects.
\end{example}

\section{Schreier invariants of groupoid extensions} \label{int th}

In this section we see coming together the two threads running parallel in the previous one,
namely that groupoid extensions are fibrations and that cocycles~/ cocycle conditions and
cohomologous cocycles are respectively lax 2-functors~/ coherence conditions and lax 2-natural
transformations. It is clear that the tying knot is the twisted product construction.

Recall from Example \ref{ex lax act and fib b.o.} that for any groupoid $\g$ the twisted product
construction establishes an equivalence of categories
$$
\semi{}\g:\LaxFun(\g,\Gr)\too \bofib(\g).
$$
Since we want to fix the fibers of our fibrations, we must fix a $\obj{\g}$-indexed family of
groups, $\k=\{K_A\}_{A\in\g}$. This family immediately determines a subcategory of
$\LaxFun(\g,\Gr)$ we already encountered, namely, the category $\Laxact(\g,\k)$ of weak actions
of $\g$ on $\Aut(\k)$. Obviously, $\Laxact(\g,\k)$ is equivalent to its image in $\bofib(\g)$ by
the above equivalence. We arrive in this way at the concept of the \emph{category of extensions
of $\g$ by $\k$}, denoted $\Ext(\g,\k)$. The problem of characterizing this category of
extensions has the following answer:

\begin{definition}\label{ext} An extension of a groupoid \g\ by a
$\obj{\g}$-indexed family of groups $\k$, is a fibration bijective on objects $P:\e\to\g$ such
that for any object $A\in \g$ the fiber of $P$ at $A$ is the group $K_{A}$. If we consider \k\ as
a totally disconnected groupoid, an extension of \g\ by \k\ is just a short exact sequence of
groupoids
$$
1\to \k\to \e\xrightarrow{P} \g\to 1,
$$
that is, \k\ is the kernel of $P$ and \g\ is the quotient groupoid of $\e$ by the normal
subgroupoid \k. The category $\Ext(\g,\k)$ is the subcategory of $\bofib(\g)$ with objects the
extensions of \g\ by \k\ and morphisms those morphisms of fibrations above \g\ which induce the
identity on fibers. A morphism of $\Ext(\g,\k)$ is then a diagram of short exact sequences of
groupoids
$$
\xymatrix@=1.5pc{1\ar[r] & \k\ar@{=}[d] \ar[r] & \e\ar[r] \ar[d] & \g\ar[r]\ar@{=}[d] &
{1\phantom{\,.}}\\ 1\ar[r] & \k \ar[r] & \e' \ar[r] & \g\ar[r] & {1\,.}}
$$
\end{definition}

Let us note that any morphism of extensions is an isomorphism and therefore $\Ext(\g,\k)$ is a
groupoid.

Using the twisted product construction it is easy to see that any weak action gives rise to an
extension. Conversely, by taking fibers, any extension gives rise to a weak action, and we have:

\begin{theorem} \label{laxact y ext} For any groupoid $\g$ and any
$\obj\g$-indexed family of groups $\k$, there is an equivalence of categories
$$
\semi {} \g:\Laxact(\g,\k)\too\Ext(\g,\k)
$$
which is the restriction to $\Laxact(\g,\k)$ of the twisted product equivalence
$\semi{}\g:\LaxFun(\g,\Gr)\to \bofib(\g)$.
\end{theorem}

According to this theorem, and with the classification problem of groupoid extensions in mind, it
makes sense to define groupoid cohomology as:

\begin{definition} \label{h2} Given a groupoid $\g$ and a family of
groups $\k$, indexed by the objects of $\g$, we define the 2-dimensional cohomology of $\g$ with
coefficients in the 2-groupoid $\Aut(\k)$ as the set of connected components of the category of
weak actions of $\g$ on $\k$
$$
H^2_{\abs{\Aut(\k)}}(\g, \k) = H^1\big(\g,\Aut(\k)\big):= \Laxact[\g,\k],
$$
We use square brackets to denote connected components. It is clear that this cohomology is not
functorial on $\k$, it is functorial in the second variable at the level of morphisms of
2-groupoids.
\end{definition} \vskip0pt

A weak action of \g\ on \k\ will be also called a \emph{2-cocycle of $\g$ with coefficients in
$\Aut(\k)$}. Then, Proposition \ref{cocycle} can be seen as a definition (or parameterization) of
2-cocycles. Two 2-cocycles will be called cohomologous if there is a morphism between them.
Example \ref{w.act map script G to a fam of gps} shows that two 2-cocycles $F_1$ and $F_2$ are
cohomologous if and only if there is $\boldsymbol{\tau}\in \mathit{C}^1(\g,\k)$ such that $F_2=
\boldsymbol{\tau} \nabla F_1$.

On the other hand, Proposition \ref{cocycle} and Example \ref{arrows between weak actions} prove
that for groups $G$ and $K$ the non abelian cohomology $H^2_{\aut(K)}(G,K)$ defined above
coincides with Dedecker \emph{\'{e}paisse} 2-cohomology of $G$ with coefficients on the
$\aut(K)$-crossed module $K$. \vskip0pt After this definition, Theorem \ref{laxact y ext} has the
following immediate corollary:

\begin{corollary}[Interpretation of groupoid cohomology] \label{inter
th} For any groupoid $\g$ and any $\obj\g$-indexed family of groups $\k$, there is a bijection
$$
H^2_{{\abs{\Aut(\k)}}}\big(\g,\k\big) \cong \Ext[\g,\k]
$$
between the 2-cocycles of $\g$ with coefficients in the 2-groupoid (crossed module) $\Aut(\k)$
and the connected components of the category $\Ext(\g,\k)$.
\end{corollary}

\section{Representation of groupoid cohomology}

The objective of this section is to give a representation theorem for the 2-dimensional non
abelian cohomology of groupoids defined in terms of homotopy classes of simplicial maps. This
theorem can be used, for example, to classify homotopy classes continuous maps from a 1-type to a
2-type in terms of this cohomology.

The key to the representation theorem (below, Theorem~\ref{repr th}) will be the construction of
nerves. We will use the \emph{geometric nerve} of a 2-category given in \cite{BuFaBl2004}. Let us
briefly recall its definition.

We denote $\sdelta$ the simplicial category, whose objects are the finite non-empty linear
orders, $\s1=[0] =\{0\}, \s2=[1] = \{0\leq1\},\dots$ and whose arrows are functors or
\emph{monotonic} maps. If $\a$ is a category and there is a functor $\embed:\sdelta\to\a$ so that
$\sdelta$ can be regarded as embedded into $\a$, one can define the $n$-simplices of an object
$A$ of $\a$ as the arrows from $\embed([n])$ to $A$, that is, the $n$-simplices of $A$ are the
$\embed([n])$-elements of $A$,
$$
\ner(A)_n = \a\big(\embed([n]),A\big),
$$
then the functoriality of $\embed$ provides face and degeneracy operators satisfying the
simplicial identities so that $\ner(A)$ becomes a simplicial set. In this way one obtains a
functor $\ner :\a \to \ssetsdel = \ssets $ defined on objects as
$$
\ner(A) = \a\big(\embed(-),A\big)
$$
and on arrows $f:A\to B$ in $\a$ via composition: $\ner(f) = \bar{f}$ is the simplicial map
$$
\ner(A)_n \xto{\ \bar{f}_n\ } \ner(B)_n, \qquad \bar{f}_n(\alpha) = f\circ\alpha.
$$

\begin{definition} If we regard $\sdelta$ as a full subcategory of
$\cat$, since $\cat$ is, in turn, a full subcategory of $\ldoscat$, we have a full embedding
$\embed : \sdelta \to \ldoscat$ and, by the above process, a functor
$$
\ner :\ldoscat \to \ssets
$$
which we take as the definition of the \emph{geometric nerve} of 2-categories.
\end{definition}

Then the nerve of a 2-category $\ba$ has:
\begin{myitemize}
\item The objects of $\ba$ as 0-simplices,
\item the arrows $A_0\xrightarrow{f} A_1$ of $\ba$ as 1-simplices,
with faces
$$
d_0(f)=A_1 \quad \text{and} \quad d_1(f)=A_0,
$$
\item the diagrams $\Delta = (g,h,f;\alpha)$ of the form

$$
\xymatrix{ & A_1 \ar@{}[d]|-(.55){\hphantom{\alpha}{\textstyle\Uparrow}\alpha} \ar[dr]^{g} \\ A_0
\ar[ur]^{f} \ar[rr]_{h} & & A_2 }
$$
(where $\alpha:h\to gf$ is a 2-cell in $\ba$) as 2-simplices, whose faces are the 1-simplices
opposite to the indicated vertex (so, in the above example, $d_0(\Delta) = g$, etc.).
\item the ``commutative'' tetrahedral $\Theta$ of the form \vskip- 1.5
cm
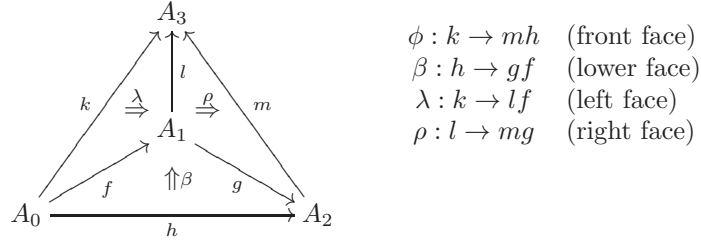
$$
\xymatrix@R-1.5pc@C-0.7pc{ & & A_3 & & \\ & & & & \\ & & & & \\ & & A_1
\ar@{}[dd]|-(.5){\hphantom{\beta}{\textstyle\Uparrow}\beta} \ar[ddrr]_-{g} \ar[uuu]_l
\ar@{}[uuull]|-(.1){\textstyle\stackrel{\lambda}{\To}}
\ar@{}[uuurr]|-(.08){\;\textstyle\stackrel{\rho}{\To}}& & \\ & & & &
\\ A_0\ar[uuuuurr]^k \ar[uurr]_-{f} \ar[rrrr]_{h} & & & &
A_2\ar[uuuuull]_m } \qquad
\begin{array}{cl} \\[1.5cm] \phi:k \to mh & \mbox{(front face)} \\
\beta:h \to gf & \mbox{(lower face)} \\ \lambda:k \to lf & \mbox{(left
face)} \\ \rho:l \to mg & \mbox{(right face)} \\
\end{array}
$$
as 3-simplices. The face operators for such tetrahedron are, as in the case of a 2-simplex, the
2-simplices opposite to the vertex indicated by the operator \big(so, for example, $d_3(\Theta) =
(m,k,h;\phi)$\big).
\end{myitemize}

For dimensions higher than 2, $\ner(\ba)$ is coskeletal. We also have that the geometric nerve
functor $\ner:\ldoscat\to\sset$ is full and faithful (see \cite[Proposition~3.3]{BuFaBl2004}) and
that a lax natural transformation between two lax functors $F,G :\ba \to \bb$ induces a homotopy
between the simplicial maps $\ner(F)$ and $\ner(G)$. Furthermore, if $\bb$ is a 2-groupoid, a
homotopy $\ner(F)\to \ner(G)$ exists if and only if there is a lax 2-natural transformation
$\alpha: G\to F$ (see \cite[Proposition~3.5]{BuFaBl2004}).

If we particularize the above results to weak actions of a groupoid \g\ on a family of groups \k,
we have that weak actions of \g\ on \k\ are in bijective correspondence with simplicial maps from
the nerve of \g\ to the nerve of $\Aut(\k)$. Moreover a morphism $\lalpha:F_1\to F_2$ between two
weak actions induces a homotopy $\bh=\ner(\lalpha):\ner(F)\to \ner(G)$ with the property that
$h_{00}:\Obj(\g)\to \Arr(|\Aut(\k)|)$ takes any object of \g\ to an identity in $|\Aut(\k)|$. We
call a \emph{normalized homotopy} to any homotopy satisfying this property, and we have

\begin{theorem}[Representation by simplicial maps] \label{repr th}
Given a groupoid $\g$ and any family $\k$ of groups indexed by the objects of $\g$, there is a
bijection
$$
H^2_{\abs{\Aut(\k)}}\big(\g,\k \big) \cong \sset_*\big[\ner(\g),\ner\big(\Aut(\k)\big)\big],
$$
where the star and the square brackets mean normalized homotopy classes of simplicial maps.
\end{theorem}

\begin{tabular}{ccc}
   V. Blanco and M. Bullejos, &  & E. Faro  \\
  Department of Algebra &  & Department of Appl. Mathematics  \\
  University of Granada &  & University of Vigo \\
  36207 Granada, Spain &  & 36207 Vigo, Spain \\
  vblanco@ugr.es, bullejos@ugr.es &  & efaro@dma.uvigo.es \\
\end{tabular}

\end{document}